\documentclass[10pt]{amsart}
\usepackage[utf8]{inputenc}
\usepackage{url,  mathrsfs}
\usepackage{hyperref}
\usepackage{amsmath}
\usepackage{graphicx, 
epstopdf}
\usepackage{subfig}
\usepackage{color}
\usepackage{bbm}
\usepackage{subfloat}
\usepackage[draft]{fixme}
\usepackage{bm} 
\usepackage{bbm} 
\usepackage{epigraph}
\usepackage{endnotes}
\usepackage{ifpdf}
\usepackage{array}

\usepackage{multicol}
\usepackage{multirow}
\usepackage{graphicx}
\usepackage{amssymb}

\theoremstyle{plain}
\newtheorem{Thm}{Theorem}
\newtheorem*{Conjecture*}{Conjecture}
\newtheorem{remark}[Thm]{Remark}

\newcommand\ul{\underline}
\newcommand{\adj}{\operatorname{adj}}

\newcommand{\R}{{\mathbb R}}

\title[The specialized V\'amos polynomial]{A note on the hyperbolicity cone of the specialized V\'amos polynomial\thanks{The author was supported by the Studienstiftung des deutschen Volkes.}
}

\begin{document}

\author{Mario Kummer}
\address{Department of  Mathematics and Statistics\\ 
University of Konstanz\\
Konstanz, Germany}
\email{mario.kummer@uni-konstanz.de}

\begin{abstract}
The specialized V\'amos polynomial is a hyperbolic polynomial of degree four in four variables with the property that none
  of its powers admits a definite determinantal representation.
  We will use a heuristical method to prove that its hyperbolicity cone
  is a spectrahedron.
\end{abstract}

\maketitle

\section{Introduction}\label{sec:intro}
 A homogeneous polynomial $h \in \R[\underline{x}]=\R[x_1,\ldots,x_n]$ is said to be hyperbolic with respect to
 $e \in \R^n$ if $h$ does not vanish in $e$ and if for every 
 $v \in \R^n$, the univariate polynomial $h(t e +v) \in \R[t]$ has only real roots.
 Originally, hyperbolic polynomials were studied in the context of partial differential equations, see for example \cite{Ga51} and \cite{HorII}.
 But recently, interest also arose in the areas of combinatorics (e.g. \cite{COS04}) and optimization (e.g. \cite{Gu97,Re06}).
 The hyperbolicity cone $\Lambda_{+}(h,e)$ of $h$ at $e$ is the set of all $v \in \R^n$ such that no zero of $h(t e +v)$ is strictly positive.
 It was shown by G{\aa}rding \cite{Gar} that $\Lambda_{+}(h,e)$ is the closure of the connected component of $\R^n \smallsetminus \{h(x)=0\}$
 which contains $e$ and that $h$ is hyperbolic with respect to every point in the interior of $\Lambda_{+}(h,e)$.
 Hyperbolicity cones are semi-algebraic convex cones, as shown for example in \cite{Gar}.
 On the other hand, a spectrahedral cone is a set 
 defined by a homogeneous linear matrix inequality, i.e. sets of the form 
 \[\{v \in \R^n: \,\, A(v)=v_1 A_1 + \ldots + v_n A_n \succeq 0\},\]
 where $A_1, \ldots, A_n$ are symmetric matrices with real entries.
 It is not hard to check that every
 spectrahedral cone is the hyperbolicity cone of an appropriate hyperbolic polynomial.
 There has been ample interest, if the converse is also true: 
 \begin{Conjecture*}
 Every hyperbolicity cone is a spectrahedral cone.
 \end{Conjecture*}
 This conjecture is commonly referred to as the Generalized Lax Conjecture. 
   Hyperbolicity cones are the feasible sets of hyperbolic programming and spectrahedral cones are the feasible sets of semidefinite programming.
   Thus the Generalized Lax Conjecture asserts that  hyperbolic programming and semidefinite programming have the same feasible sets.
 An equivalent formulation of this conjecture is:
 \begin{Conjecture*}
  Let $h \in \R[\underline{x}]$ be hyperbolic with respect to $e \in \R^n$. Then there is a hyperbolic polynomial $q \in \R[\underline{x}]$, such
  that $\Lambda_{+}(h,e) \subseteq \Lambda_{+}(q,e)$ and such that $q \cdot h$ has a definite determinantal representation, i.e.
  \[q \cdot h = \det A(\underline{x})=\det(x_1 A_1 +\ldots +x_n A_n),\] with $A_1, \ldots , A_n$ symmetric matrices with real entries and
  $A(e) \succ 0$.
 \end{Conjecture*}
 It is easy to see that this formulation implies the other one, because the hyperbolicity cone of a product of hyperbolic polynomials
 is just the intersection of the hyperbolicity cones of the factors. For the other direction we refer to \cite{HelVin07}.
 In the case $n=3$ the conjecture is true  \cite{HelVin07} and also in the case when $h$ is an elementary symmetric polynomial \cite{Bran}.
 In general, if the projective variety defined by $h$ is smooth, one can always find a hyperbolic polynomial $q$
 such that $q \cdot h$ has a definite determinantal representation \cite{Kum} (without control on the hyperbolicity cone of $q$).
 
 We will explain a heuristical method to find a representation of a given hyperbolicity cone as a spectrahedral cone in Section \ref{sec:allgkonstr}.
 In Section \ref{sec:vamos}, we will apply this method to the specialized V\'amos polynomial $h_4$, which is a hyperbolic polynomial of degree four in four
 variables, with the property that no power $h_4^N$ has a definite determinantal representation.
 Note that only very few polynomials \cite{Bra11,Burt14} are known to have this property and all of these come from the V\'amos matroid.
 The polynomial $h_4$ and its relatives frequently serve as counterexamples for questions concerning hyperbolic polynomials \cite{Bra11,Us,NPT}
 and thus it is natural to check whether it might also be a counterexample to the Generalized Lax Conjecture. The main result of this note is that
 indeed it is not such a counterexample.

\section{The construction method}\label{sec:allgkonstr}
 We denote by $\R[\underline{x}]_d$ the set of all homogeneous forms in $\R[\underline{x}]$ of degree $d$.
 Let $h \in \R[\underline{x}]_d$ be irreducible and hyperbolic with respect to $e \in \R^n$.
 Let $\underline{f} \in \R[\underline{x}]_{d'}^m$, where $d' \geq d-1$.
 Further, we want $\underline{f}$ to have at least one entry that is not divisible by $h$.
 Let $A_1, \ldots, A_n$
 be symmetric matrices with real entries of size $m$ and let $\underline{g} \in \R[\underline{x}]_{d'-d+1}^m$.
 Consider the conditions
 \begin{equation}\label{eq:sdp}
  A(\underline{x}) \cdot \underline{f} = h \cdot \underline{g}, \,\, A(e) \succ 0.
  \end{equation}
  If $h$ and $\ul{f}$ are fixed, this corresponds to a semidefinite program in the entries of the $A_i$ and
  the coefficients of the entries of $\ul{g}$. 
 If (\ref{eq:sdp}) has a solution, then there is a hyperbolic polynomial $q \in \R[\ul{x}]_{m-d}$ such that $q \cdot h$ has
 a definite determinantal representation, namely $\det A(\ul{x})=q \cdot h$.
 If we are lucky, then there is a solution such that even $\Lambda_{+}(h,e) \subseteq \Lambda_{+}(q,e)$ holds.
 The proof of the main result of \cite{Kum} shows that (\ref{eq:sdp}) always has a solution 
 if we choose $\ul{f}$ to be the vector of all monomials of high enough degree and
 if the hypersurface in projective space
 defined by $h$ is smooth.
 
 Note that, conversely, if there is a hyperbolic polynomial $q \in \R[\ul{x}]$ such that $q \cdot h$ has a definite determinantal representation
 $q \cdot h = \det A(\ul{x})$ with symmetric matrices $A_i$ of size $m$, then there is a $\ul{f} \in \R[\ul{x}]_{d'}^m$ for some $d' \geq d-1$ such
 that every entry of $A(\ul{x})\cdot \ul{f}$ is divisible by $h$, but not every entry of $\ul{f}$.
 Take, for example, some row of the matrix $\adj (A(\ul{x}))$ as $\ul{f}$.
 
 \section{The V\'amos polynomial}\label{sec:vamos}
 Consider the specialized V\'amos polynomial: \[h_4=x_1^2 x_2^2 + 
  4 \cdot (x_1 + x_2 + x_3 + x_4)\cdot (x_1 x_2 x_3 + x_1 x_2 x_4 + x_1 x_3 x_4 + x_2 x_3 x_4).\]
  Wagner and Wei \cite{WW} have shown that $h_4$ is hyperbolic with respect to \[e=(1,1,0,0)^{\rm T}\] and
  Br\"and\'en \cite{Bra11} has proved that no power of $h_4$ admits a definite determinantal representation.
\begin{remark}\label{rem:vamosecht}
  Note that in fact the authors of \cite{Bra11} and \cite{WW} considered the polynomial $h_{V_8}$, which is the bases generating polynomial
  of the V\'amos matroid $V_8$. This is a multiaffine polynomial in eight variables of degree four and one can obtain $h_4$ by restricting
  $h_{V_8}$ to a linear subspace. Wagner and Wei \cite{WW} proved that $h_{V_8}$ is hyperbolic and Br\"and\'en proved that no power of $h_{V_8}$
  admits a definite determinantal representation.
  It is easy to see, that both properties spread to $h_4$.
  All polynomials that are known to have these properties are constructed in some way from the V\'amos matroid, see also \cite{Burt14}.
\end{remark}
 We will show that the hyperbolicity cone of $h_4$ is nevertheless a spectrahedral cone.
 In order to apply the construction of the previous section, we choose $d'=3$, $m=7$ and
 \[\ul{f}=\left(
\begin{array}{c}
 x_1^2 x_2 \\
 x_1 x_2^2+x_1 x_2 x_3 \\
 x_1 x_2 x_4 \\
 x_1 x_3 x_4+x_2 x_3 x_4+x_3^2 x_4+x_3 x_4^2 \\
 x_2^2 x_4+x_2 x_3 x_4+x_2 x_4^2 \\
 x_1^2 x_4+x_1 x_3 x_4+x_1 x_4^2 \\
 -x_1 x_2^2+x_1^2 x_3+x_1 x_3^2+x_1 x_3 x_4
\end{array}
\right).\]
\begin{remark}
 All entries of $\ul{f}$ vanish on the real singularities of the hypersurface defined by $h_4$.
\end{remark}
 Solving the semidefinite program (\ref{eq:sdp}) (for example using \cite{yalmip}) yields the following four symmetric $7 \times 7$ matrices:
 \[A_1=\left(
\begin{array}{ccccccc}
 0 & 0 & 0 & 0 & 0 & 0 & 0 \\
 0 & 1 & 2 & 0 & 1 & 0 & 0 \\
 0 & 2 & 6 & 0 & 4 & 0 & 0 \\
 0 & 0 & 0 & 4 & 0 & 0 & 0 \\
 0 & 1 & 4 & 0 & 4 & 0 & 0 \\
 0 & 0 & 0 & 0 & 0 & 0 & 0 \\
 0 & 0 & 0 & 0 & 0 & 0 & 0
\end{array}
\right)
  ,\,
  A_2=\left(
\begin{array}{ccccccc}
 1 & 1 & 2 & 1 & 0 & 2 & 1 \\
 1 & 4 & 4 & 4 & 0 & 4 & 4 \\
 2 & 4 & 9 & 4 & 0 & 8 & 4 \\
 1 & 4 & 4 & 4 & 0 & 4 & 4 \\
 0 & 0 & 0 & 0 & 0 & 0 & 0 \\
 2 & 4 & 8 & 4 & 0 & 8 & 4 \\
 1 & 4 & 4 & 4 & 0 & 4 & 4
\end{array}
\right),
 \]
       \[
  A_3=\left(
\begin{array}{ccccccc}
 3 & 3 & 10 & 0 & 3 & 4 & 0 \\
 3 & 4 & 12 & 0 & 4 & 4 & 0 \\
 10 & 12 & 41 & 0 & 12 & 16 & 0 \\
 0 & 0 & 0 & 0 & 0 & 0 & 0 \\
 3 & 4 & 12 & 0 & 4 & 4 & 0 \\
 4 & 4 & 16 & 0 & 4 & 8 & 0 \\
 0 & 0 & 0 & 0 & 0 & 0 & 0
\end{array}
\right),
  A_4=\left(
\begin{array}{ccccccc}
 2 & 2 & 2 & 0 & 0 & 0 & 0 \\
 2 & 7 & 3 & 0 & 0 & 0 & 4 \\
 2 & 3 & 4 & 0 & 0 & 0 & 0 \\
 0 & 0 & 0 & 0 & 0 & 0 & 0 \\
 0 & 0 & 0 & 0 & 0 & 0 & 0 \\
 0 & 0 & 0 & 0 & 0 & 0 & 0 \\
 0 & 4 & 0 & 0 & 0 & 0 & 4
\end{array}
\right)
  .
 \]
 The matrix $A_1+A_2$ is positive definite, all eigenvalues are bigger than $0,1$.
 We will show that the hyperbolicity cone $\Lambda_{+}(h_4,e)$ has the following representation:
 \[\Lambda_{+}(h_4,e)=\{x \in \R^4: \,\,\, x_1 A_1+x_2 A_2+x_3 A_3+x_4 A_4 \succeq 0\}.\]
 Using a computer algebra system, one can verify that
 \[\det(x_1 A_1+x_2 A_2+x_3 A_3+x_4 A_4)=32 \cdot l \cdot q \cdot h_4 ,\]where
 $l=2x_1+3x_2+3x_3+4x_4$ and $q=x_1 x_2+ x_1 x_3+2x_1 x_4+ x_2 x_4 + x_3 x_4$.
 It remains to verify that the hyperbolicity cone of $l\cdot q$ contains the hyperbolicity cone of $h_4$.
 We have that $A(\underline{x}) \cdot \underline{f} = h_4 \cdot \underline{g}$ where $\underline{g}=(1,2,4,1,1,2,1)^{\rm T}$.
 Consider the polynomial 
 $p=\underline{g}^{\rm T} \underline{f} .$ One checks that \[p= (x_1+x_2+x_3+x_4) \cdot q.\]
 Applying a simple calculation, we see that
   \[\textrm{D}_e h_4 \cdot p- h_4 \cdot \textrm{D}_e p=\sum_{i=1}^2(\frac{\partial h_4}{\partial x_i} \cdot p - h_4 \cdot \frac{\partial p}{\partial x_i} )= \ul{f}^{\rm T} \cdot( A_1+A_2) \cdot \ul{f}.\]
 In particular, since $A_1+A_2$ is positive definite, this mixed derivative is nonnegative on $\R^4$.  
 It follows from \cite[Theorem 2.1]{Us} that $p$ is hyperbolic with respect to $e$  and that its hyperbolicity cone contains the hyperbolicity cone of $h_4$.
 In particular, we have $\Lambda_{+}(h_4,e) \subseteq \Lambda_{+}(p,e) \subseteq \Lambda_{+}(q,e)$.
 Since $q$ is quadradic and since $q(2,0,0,1)>0$ we have that $q$ is hyperbolic with respect to $(2,0,0,1)^{\rm T}$ and thus the derivative
 of $q$ in that direction, which is $l$, is nonnegative on $\Lambda_{+}(q,e)$. Therefore, we indeed have
 \[
  \Lambda_{+}(h_4,e) \subseteq \Lambda_{+}(q,e) = \Lambda_{+}(l \cdot q,e).
 \]
This shows that the hyperbolicity cone $\Lambda_{+}(h_4,e)$ is a spectrahedral cone.

\bibliographystyle{plain}
\bibliography{vamos}

\end{document}